\documentclass[12pt, a4paper, oneside]{amsart}
\usepackage[a4paper, left=2cm, right=2cm, top=2cm, bottom=2cm]{geometry}
\usepackage{a4wide}
\usepackage[english]{babel}
\usepackage[T2A]{fontenc}
\usepackage[utf8]{inputenc} 
\usepackage{amsfonts}
\usepackage{amssymb, amsthm, amscd}
\usepackage{amsmath}
\usepackage{mathtools}
\usepackage{needspace}
\usepackage{etoolbox}
\usepackage{lipsum}
\usepackage{comment}
\usepackage{cmap}
\usepackage[pdftex]{graphicx}
\usepackage[unicode]{hyperref}
\usepackage[matrix,arrow,curve]{xy}
\usepackage[usenames,dvipsnames]{xcolor}
\usepackage{colortbl}
\usepackage{textcomp}
\usepackage{cite}
\usepackage{euscript}

\DeclareFontFamily{U}{mathx}{}
\DeclareFontShape{U}{mathx}{m}{n}{<-> mathx10}{}
\DeclareSymbolFont{mathx}{U}{mathx}{m}{n}
\DeclareMathAccent{\widehat}{0}{mathx}{"70}
\DeclareMathAccent{\widecheck}{0}{mathx}{"71}

\makeatletter
\def\@settitle{\begin{center}
		\baselineskip14\p@\relax
		\bfseries
		\LARGE
		\@title
	\end{center}
}
\makeatother

\newtheorem{theorem}{Theorem}[section]

\newtheorem{lemma}[theorem]{Lemma}
\newtheorem{corollary}[theorem]{Corollary}

\newtheorem{definition}[theorem]{Definition}
\newtheorem{mainth}{Theorem}

\theoremstyle{definition}

\renewcommand{\geq}{\geqslant}
\renewcommand{\leq}{\leqslant}
\newcommand{\R}{\mathbb{R}}
\newcommand{\N}{\mathbb{N}}
\newcommand{\Z}{\mathbb{Z}}

\newcommand{\G}{\mathcal{G}}
\newcommand{\T}{\mathbb{T}}
\newcommand{\LL}{\mathcal{L}}
\newcommand{\eps}{\varepsilon}

\DeclareMathOperator{\Ker}{Ker}
\DeclareMathOperator{\dist}{dist}

\DeclareMathOperator{\ch}{ch}

\newenvironment{enumerate*}%
  {\begin{enumerate}%
    \setlength{\itemsep}{1pt}%
    \setlength{\parskip}{1pt}}%
  {\end{enumerate}}

\makeatletter
\newcommand*{\@old@slash}{}\let\@old@slash\slash
\def\slash{\relax\ifmmode\delimiter"502F30E\mathopen{}\else\@old@slash\fi}
\makeatother

\title{Irregular Gabor Frames for the Poisson Kernel}
\author{Yurii Belov and Andrei V. Semenov}

\thanks{The work was supported by RSF grant 26-11-00157}

\address{Andrei V. Semenov
Saint Petersburg State University\\
Russia, 199178, Saint Petersburg, 14 line of Vasilievsky island, 29}

\email{asemenov.spb.56@gmail.com}

\address{Yurii Belov
Saint Petersburg State University\\
Russia, 199178, Saint Petersburg, 14 line of Vasilievsky island, 29}

\email{j\_b\_juri\_belov@mail.ru }

\begin{document}
\maketitle

\begin{abstract}
Consider the semi-regular system of time-frequency shifts of the Poisson kernel $g(t) = {{1} \over {t^2+w^2}}$
$$\G(g, M \times \alpha \Z) = \left\{ e^{2\pi i \lambda t } g(t-n) \colon (\lambda, n) \in M \times \alpha\Z \right\}$$
 where $M \subset \R$. We obtain {\it a criterion} for $\G(g, M \times \alpha \Z)$ to form a frame in $L^2(\R)$ for an {\it arbitrary} set $M$ of finite upper density. It turned out that the lower density alone does not determine the frame property, even when $M$ is uniformly discrete.
\end{abstract}

\section{Introduction}
For a function $g \in L^2(\R)$ and a point $(\lambda, \mu) \in \R^2$, we define the time-frequency shift $\pi_{\lambda, \mu} g$ by 
$$\pi_{\lambda, \mu} g (t)= e^{2\pi i \lambda t} g(t - \mu).$$
Now for a discrete set $\LL \subset \R^2$ we consider a {\it Gabor system} $\G(g, \LL)$:
$$\G(g, \LL) := \{\pi_{\lambda, \mu} g \mid (\lambda, \mu) \in \LL\}, \quad \text{where } g\in L^2(\R).$$
We say that $\G(g,\LL)$ is a {\it Gabor frame} if for some constants $A, B >0$ one has
\begin{equation}
A\|f\|^2_2\leq \sum_{(\lambda, \mu) \in \LL}|\langle f, \pi_{\lambda, \mu}g \rangle|^2\leq B\|f\|^2_2,  \text{ for any } f\in L^2(\mathbb{R}),
\label{frameineq}
\end{equation}
where $\langle \cdot, \cdot \rangle$ denotes the standard inner product on $L^2(\R)$. Inequality (\ref{frameineq}) is referred to as {\it the frame inequality}. \par 
The most natural choice of $\LL$ is a rectangular lattice $\alpha \Z \times \beta\Z$. The problem of finding the frame set, i.e., all pairs $(\alpha, \beta)$ such that $\G(g,\alpha \Z \times \beta \Z)$ is a frame in $L^2(\R)$ is crucial in time-frequency analysis.\par
We study the Poisson kernel $g(t) = {{1} \over {t^2 + w^2}}$. In the classical case of $\LL=\alpha \Z \times \beta \Z$ the frame criterion is known:
$\G(g,\alpha\Z\times\beta\Z)$ forms a frame in $L^2(\R)$ if and only if $\alpha \beta < 1$ (see \cite{Gro2}). This paper is devoted to the semi-regular case $\LL = M \times \beta \Z$.

\subsection{State of the Art}

The classical case $\LL = \alpha \Z \times \beta \Z$ for $\alpha, \beta >0$ is the most
studied. For a detailed overview we refer to the monograph \cite{Gro} and to the surveys
\cite{Heil, Mystery, GK}. A complete description of the set of parameters $(\alpha, \beta) \in \R^2$ such that $\G(g, \alpha \Z \times \beta \Z)$ forms a frame in $L^2(\R)$ is known for totally positive functions of finite type at least two (see \cite{Gro2}) and for totally positive functions of Gaussian type (see \cite{Gro1}). In both cases, $\G(g,\alpha\Z\times\beta\Z)$ is a frame in $L^2(\R)$ if and only if $\alpha\beta<1$. Earlier complete characterizations were available for Gaussian windows (see \cite{L,S,SW}), exponential windows (see \cite{Jans2,Jans3}), and the hyperbolic secant (see \cite{JansStr}). Also, full characterizations were obtained for rational functions of Herglotz type (see \cite{BKL1}), shifted sinc-functions (see \cite{BAVS}),  hyperbolic secant type functions (see \cite{BB}) and nonnegative functions supported on a half-line and satisfying the uniform decay condition of \cite{BK}. \par
Very recently, J. de Dios Pont, K.~Gr\"ochenig, L.~Liehr, I. Shafkulovska and M. Taylor proved that for every continuous totally positive function $g\in L^1(\R)$ the system $\G(g,\alpha\Z\times\beta\Z)$ is a frame in $L^2(\R)$ if and only if $\alpha\beta<1$. This settles the rectangular-lattice frame set conjecture for
continuous integrable totally positive functions, thus concluding a research program that had lasted for many years, see \cite{PGLST}.\par 

When moving to non-classical index sets $\LL \subset \R^2$, density remains a fundamental necessary constraint (see \cite{RS}). It is known that if $\LL \subset \R^2$ is an arbitrary discrete set with $D_{-}(\LL)<1$, then $\G(g,\LL)$ is not a frame for any $g\in L^2(\R)$, see \cite[Theorem 13(b)]{Heil}. For windows in Feichtinger's algebra $M^1(\R)$, the necessary inequality is strict, but it is not sufficient in general (see \cite[Theorem~1.5 and Remark~1.6(ii)]{AFK}). For a Gaussian window and a uniformly discrete set $\LL\subset\R^2$, however, $D_{-}(\LL)>1$ is both necessary and sufficient (see \cite{L,S,SW}).

For a totally positive window $h$ of Gaussian type and a uniformly discrete set $\Lambda\subset\R$, the semi-regular system $\G(h,\beta\Z\times\Lambda)$ is a frame if and only if
$D_{-}(\Lambda)>\beta$, using notation of this paper (see \cite[Theorem~6.2]{Gro1}).
Fourier transformation gives the corresponding criterion for $\G(\widehat h,\Lambda\times\beta\Z)$. Under the same separation assumption, sharp density criteria
for both semi-regular orientations are known for the hyperbolic secant-type windows of \cite[Theorem~1.1]{BB}. For the Cauchy window $(t-iw)^{-1}$, Belov, Kulikov, and
Lyubarskii obtained an exact frame criterion for direct-product index sets in terms of sampling in
a Paley--Wiener space (see \cite{BKL2}).

We give a complete geometric criterion for semi-regular Gabor frames generated by the Poisson kernel. Although this window belongs to $M^1(\R)$, the frame property of $\G(g,M\times\alpha\Z)$ is not determined by the lower density of $M$ alone, even when $M$ is uniformly discrete.
In particular, a gap of length at least $2/\alpha$ precludes the frame property, regardless of the density of $M$ outside that gap. Moreover, the answer depends on the mutual position of the gaps of $M$ modulo $1/\alpha$, as Condition (3) of Theorem \ref{th:main} shows. \par

Recently in \cite{AVS} the second author proved that there exists a universal set of the form $M \times \Z$ for a wide class of well-behaved rational functions. We develop a similar technique in a more general setting.\par

\subsection{Main result}
Consider the Poisson kernel $g(t) = {{1} \over {t^2 +w^2}}$ for some $w>0$: the simplest example of a rational function that satisfies the conditions of the Balian–Low theorem. For a system $\G(g, M \times \alpha\Z)$, we consider the problem of finding necessary and sufficient conditions for such a system to generate a frame. Note that for $M = \beta\Z$ the problem is classical and completely solved.
\begin{definition} A set $M \subset \R$ is called {\it uniformly discrete} if there exists $\delta> 0 $ such that $|x-y| > \delta$ for any distinct $x, y \in M$.
\end{definition}
 It is well-known that if $\G(g, M \times \Z)$ is a frame then there exists a uniformly discrete subset $L\subseteq M$ such that $\G(g, L \times \Z)$ is also a frame (see [\cite{GOR}, Theorem 7.1]). \par 
Recall that the {\it lower Beurling density} $D_{-}(\Lambda)$ is defined as:
$$D_{-}(\Lambda) = \liminf_{a \to \infty} \inf_{R \in \R} {{\#\{x\in \Lambda \mid x \in [R, R+a]\}} \over {a}}.$$
The upper Beurling density is defined by
\[D^+(\Lambda)=\limsup_{a\to\infty}\sup_{R\in\mathbb R}
 {{\#\{x\in \Lambda \mid x \in [R, R+a]\}} \over {a}}.
\]

Whenever a set is indexed by $\mathbb Z$ below, its points are listed in strictly increasing order. Our goal is to prove a criterion for an arbitrary set $M$ of finite upper density: this condition is already necessary for the upper bound in \eqref{frameineq}.\par
 For arbitrary $a,b\in\mathbb{R}$ let
$$[a,b)_{\alpha^{-1} \mathbb{T}}=\{x \mid x\in[0,\alpha^{-1}), \text{ for some } n\in\mathbb{Z}, a\leq x+n\alpha^{-1}<b\},$$
and define $\dist(X;Y)=\inf_{x\in X, y\in Y}|x-y|$. In particular, $[a,b)_{\alpha^{-1} \mathbb{T}}=\emptyset$ for $b\leq a$, $\dist(\emptyset; Y)=\infty$ for any set $Y$.

\begin{mainth}\label{th:main} {\it Let $M \subset \R$ be of finite upper density and let $\alpha>0$. The system $\G\left({{1} \over {t^2+w^2}}, M \times \alpha\Z\right)$ forms a frame in $L^2(\R)$ if and only if there exists $\eps>0$ and a uniformly discrete subset $N = \{n_k\}_{k \in \Z} \subseteq M$ such that
\begin{enumerate}
\item $D_{-}(N) > \alpha$,
\item $\sup_{k \in \Z} |n_{k+1} - n_k| < {{2} \over {\alpha}} $,
\item for any indices $m < k$ such that $n_{m+1} - n_m > 1/\alpha$ and $n_{k+1} - n_k > 1/\alpha$, the following intersection is empty:
$$\left[n_m  - \eps, n_{m+1}-{{1} \over {\alpha}}\right)_{\alpha^{-1} \mathbb{T}} \cap \left[n_k, n_{k+1}-{{1} \over {\alpha}} + \eps\right)_{\alpha^{-1} \mathbb{T}} \cap \left( \bigcap_{j=m+1}^{k-1} [n_j, n_{j+1})_{\alpha^{-1} \mathbb{T}} \right) = \emptyset.$$


\end{enumerate}
}
\end{mainth}

Thus neither the Beurling densities nor the gap lengths determine the frame property: the arrangement of the gaps is essential, even for finite modifications of arbitrarily
dense lattices.\par 

We want to point out the crucial non-symmetry of the problem. While the criterion for the time-frequency sets of the form $M \times \alpha \Z$ is established, the authors do not know the answer for the sets of the form $\alpha \Z \times M$. Any progress in this direction would be of great interest to us.

 \subsection*{Organization of the paper} In Section 2 we formulate and prove the main criterion, which allows us to reduce the problem to a uniform lower bound for a family of matrix operators $L_\xi$. In Section 3 we give a brief outline of the paper and construct the operator $L_{\xi}$ explicitly. In Section 4 we discuss the inner structure of the operator $L_{\xi}$ constructed in Section 3. In Section 5 we prove the necessity of Theorem \ref{th:main}. In Section 6 we prove the sufficiency of Theorem \ref{th:main}.

 \subsection*{Notations} We normalize the Fourier transform as follows: $\widehat{f}(z) =\displaystyle \int_{\R} f(t) e^{2\pi i t z} dt.$
We let $PW_a$ denote the Fourier transform of $L^2(0,a)$. By $C_0^{\infty}(\R)$ we denote the space of all infinitely differentiable functions on $\R$ with compact support. By $B_R(a)$ we denote a ball of radius $R$ centered at $a$. For $x \in \R$ denote the integer part of $x$ as $[x]$ and the fractional part of $x$ as $\{x\}$. We define $\lceil x\rceil$ to be the smallest integer not less than $x$.\par
Finally, by $U(z) \lesssim V(z)$ (equivalently $V(z) \gtrsim U(z)$), we mean that there exists a constant $C>0$ such that $U(z) \leq C V(z)$ holds for all $z$ in the set under consideration. We write $U(z) \asymp V(z)$ if both $U(z) \lesssim V(z)$ and $U(z) \gtrsim V(z)$ hold.

\section{Main criterion}
The technique developed here is similar to the one used in the proof of Theorem 2.1 in \cite{AVS}. Let $M = \{m_n \mid n \in \Z\}$ be an increasing sequence, i.e. $m_n < m_{n+1}$ for all $n \in \Z$. We define  
$$\pi_{k,n}  g(t) = \pi_{m_k,n} g(t)=e^{2\pi i m_k t} g(t-n)$$
for all $k, n \in \Z$. In what follows, we use the symbol $\pi_{k,n}$ in the sense defined above.\par

\begin{theorem}\label{main:crit} Let $M\subset \R$ be uniformly discrete. The system $\G\left({{1} \over {t^2+w^2}}; M \times \alpha \Z\right)$ forms a frame in $L^2(\R)$ if and only if
\begin{equation}
\label{maincr}
\|G\|_2^2 \asymp \sum_{k \in \Z}  \int_0^{1/\alpha} \left| G(t+m_k) \sinh (2\pi wt) - G\left(t+m_k+{{1} \over {\alpha}}\right) \sinh \left(2\pi w\left(t - {{1} \over {\alpha}}\right)\right) \right|^2 dt
\end{equation}
for any $G \in L^2(\R)$.
\end{theorem}
\begin{proof}

Without loss of generality, we may assume that $\alpha = 1$. Indeed, for $a>0$ put $U_a f(t) = a^{1/2} f(at)$. The operator $U_a$ is unitary on $L^2(\R)$ and
$$U_a \pi_{\lambda, \mu} g = a^{1/2} \pi_{a \lambda, \, \mu / a} \left( g(a \, \cdot) \right), \qquad (\lambda, \mu) \in \R^2,$$
so that $\G(g; M \times \alpha \Z)$ forms a frame in $L^2(\R)$ if and only if $\G(g(\alpha \, \cdot); \alpha M \times \Z)$ does. 

Since $g(\alpha t)=\alpha^{-2}\frac{1}{t^2+(w/\alpha)^2}$, this rescaling replaces $M$ and $w$ by $\alpha M$ and $w/\alpha$, respectively. We therefore first prove the
equivalence for $\alpha=1$. It suffices to show that
$$\| f\|_2^2 \asymp \sum_{k \in \Z} \sum_{n \in \Z}  |\langle \pi_{k,n} g, f \rangle |^2 \text{ for all } f \in C_0^{\infty}(\R).$$
Let $\mathcal C= \left\{ c\in\ell^2(\Z\times\Z) \mid \operatorname{supp}c\text{ is finite} \right\}$. Fix $f\in C_0^\infty(\R)$. By duality, we have
$$ \left( \sum_{k \in \Z} \sum_{n \in \Z}  |\langle \pi_{k,n} g, f \rangle |^2 \right)^{1/2} = \sup\limits_{c \in \mathcal{C} \colon \|c\| = 1} \left| \sum_{k \in \Z} \sum_{n \in \Z} c_{k,n} \langle \pi_{k,n} g, f \rangle \right|.$$

Rewriting the inner sum, we obtain
$$\sup\limits_{c \in \mathcal C \colon \|c\| = 1}\left|  \sum_{k \in \Z} \sum_{n \in \Z} c_{k,n} \langle \pi_{k,n} g, f \rangle \right| = \sup\limits_{c \in \mathcal C \colon \|c\| = 1} \left| \sum_{k \in \Z} \sum_{n \in \Z} c_{k,n} \int_{\R} e^{2\pi i m_k t} {{\overline{f(t)}} \over {(t-n)^2+w^2}}  dt \right| \asymp$$
$$\asymp \sup\limits_{c \in\mathcal C \colon \|c\| = 1}\bigg|  \int_{\R} {{\overline{f(t)} } \over {(1-e^{2\pi i (t-iw)})(1-e^{2\pi i (t+iw)})}} \sum_{k \in \Z} e^{2\pi i m_k t} \times$$
$$\times (1-e^{2\pi i (t-iw)})(1-e^{2\pi i (t+iw)}) \sum_{n \in \Z} c_{k,n} \left({{1} \over {t-(n-iw)}} - {{1} \over {t-(n+iw)}} \right) dt \bigg|.$$
We define a function $F \in L^2(\R)$ by the rule
$$\overline{F(t)} = {{\overline{f(t)} } \over {(1-e^{2\pi i (t-iw)})(1-e^{2\pi i (t+iw)})}}.$$
Put $G=\widehat{\overline F}$. Note that $\|G\|_2 \asymp \|F\|_2 \asymp \|f\|_2$ since the denominator of $F$ is bounded away from zero and infinity. For each $k \in \Z$, let
$$g_k (z)= (1-e^{2\pi i z}) \sum_{n \in \Z}  {{c_{k,n}} \over {z-n}}.$$
Since $\{c_{k,n}\}_{n \in \Z} \in \ell^2(\Z)$, it follows that $g_k \in PW_1$ for each fixed $k \in \Z$. This yields
$$\left( \sum_{k, n \in \Z} | \langle \pi_{k,n} g, f\rangle |^2 \right)^{1/2} \asymp \sup\limits_{c \in \mathcal C \colon \|c\| = 1} \bigg| \int_{\R} \overline{F(t)}  \times $$
$$\times \sum_{k \in \Z} e^{2\pi i m_k t} \left( (1-e^{2\pi i (t+iw)})g_k(t-iw) - (1-e^{2\pi i (t-iw)})g_k(t+iw) \right) dt \bigg|.$$
For $k\in\mathbb Z$ and $0<t<1$, put
\[Q_k(G;t)=G(t+m_k)\sinh(2\pi wt)-G(t+m_k+1)\sinh\bigl(2\pi w(t-1)\bigr).
\]
We regard $Q(G)=\{Q_k(G;\cdot)\}_{k\in\mathbb Z}$ as an
element of $H = \ell^2\bigl(\mathbb Z;L^2(0,1)\bigr)$. The uniform discreteness of $M$ implies that the intervals $[m_k,m_k+2]$ have bounded overlap. Consequently,
\[
\begin{aligned}
\|Q(G)\|_{H}^2&=
\sum_{k\in\mathbb Z}\int_0^1|Q_k(G;t)|^2\,dt\leq 2\sinh^2(2\pi w)\sum_{k\in\mathbb Z}\int_{m_k}^{m_k+2}|G(u)|^2\,du \lesssim \|G\|_2^2.
\end{aligned}
\]

For a finitely supported coefficient sequences we obtain
\[
\check g_k(t)=-2\pi i\sum_{n\in\mathbb Z}c_{k,n}e^{-2\pi int} \quad \text{for every } t \in (0,1).
\]
 Hence $\sum_{k\in\mathbb Z}\|\check g_k\|_{L^2(0,1)}^2=4\pi^2\sum_{k,n\in\mathbb Z}|c_{k,n}|^2$. Now 
\[
\int_{\mathbb R}\overline{F(u)}e^{2\pi i m_k u}g_k(u\pm iw)\,du=\int_0^1G(t+m_k)\check g_k(t)e^{\mp2\pi wt}\,dt.
\]
Now we obtain
\[
\begin{aligned} \left( \sum_{k,n\in\Z} |\langle\pi_{k,n}g,f\rangle|^2 \right)^{1/2}
\asymp \sup_{\substack{c\in\mathcal C\\\|c\|_{\ell^2}=1}} \left| \sum_{k\in\Z}\int_0^1 \check g_k(t)Q_k(G;t)\,dt \right| \asymp \|Q(G)\|_H.
\end{aligned}
\]
The last comparison follows by duality, since the sequences $(\check g_k)_k$ corresponding to unit vectors in $\mathcal C$ are dense in the sphere of radius $2\pi$ in $H$. Since $f\mapsto G$ is a bounded invertible antilinear map on $L^2(\R)$ and $\|G\|_2\asymp\|f\|_2$, this completes the proof.
\end{proof}

\section{Outline of the proof of Theorem \ref{th:main}}
One can assume that $\alpha = 1$ for the reasons given in the proof of Theorem \ref{main:crit}. Note that this reduction replaces $w$ by $w/\alpha$, so that $w>0$ below is an arbitrary fixed parameter, and all the constants may depend on it. Under the same rescaling, a pair $(N,\eps)$ satisfying conditions (1)--(3) for $M$ and $\alpha$ corresponds to the pair $(\alpha N, \alpha\eps)$ satisfying conditions (1)--(3) for $\alpha M$ and $\alpha = 1$. \par
From now on, we consider only the time-frequency set $M \times \Z$. The upper estimate in  (\ref{maincr}) follows from the bounded overlap
of the intervals $[m_k,m_k+2]$ and the boundedness of the coefficients on $[0,1]$. Thus, we only need to verify the lower estimate.

\medskip

Consider the unitary map $U: L^2(\mathbb R) \longrightarrow L^2\bigl([0,1);\ell^2(\mathbb Z)\bigr)$ defined by the formula
\[ (U G)(\xi) = \{(-1)^jG(\xi+j)\}_{j\in\mathbb Z}. \]
The sequence is defined for almost every $\xi$, and the map is understood on equivalence classes of functions. It is the usual decomposition into translates of $[0,1)$, composed with the diagonal unitary map
$\{x_j\}\mapsto\{(-1)^jx_j\}$. In particular,
\[
\|G\|_{2}^2 = \int_0^1 \|(U G)(\xi)\|_{\ell^2(\mathbb Z)}^2\,d\xi.
\]

For $\xi\in[0,1)$ and $k\in\mathbb Z$ put
\[
t_k=\{\xi-m_k\} \quad \text{and} \quad b_k=t_k+m_k-\xi=\lceil m_k-\xi\rceil.
\]
Define the operator $L_\xi:\ell^2(\mathbb Z)\to\ell^2(\mathbb Z)$ by the rule
\[ (L_\xi x)_k = \sinh(2\pi wt_k)x_{b_k} + \sinh\bigl(2\pi w(t_k-1)\bigr)x_{b_k+1}.
\]
Uniform discreteness now gives $\sup_{\xi\in[0,1)}\|L_\xi\|<\infty$. We assign the column pair $(b_k,b_k+1)$ to the $k$-th row, even if one of the corresponding coefficients vanishes. With $Q_k$ as defined in the proof of Theorem \ref{main:crit}, we have
\[
\bigl(L_\xi(U G)(\xi)\bigr)_k = (-1)^{b_k}Q_k(G;t_k)
\]
for almost every $\xi$. Changing variables $\xi=\{t+m_k\}$ separately in each
summand, we obtain
\[
\sum_{k\in\mathbb Z}\int_0^1|Q_k(G;t)|^2 dt =\int_0^1 \|L_\xi(U G)(\xi)\|_{\ell^2(\mathbb Z)}^2 d\xi.
\]
It now follows from Theorem \ref{main:crit} that
$\G(g;M\times\mathbb Z)$ is a frame if and only if there exists $c>0$ such that
\begin{equation}\label{eq:opercrit}
\|L_\xi x\|_{\ell^2(\mathbb Z)}^2 \ge c\|x\|_{\ell^2(\mathbb Z)}^2
\quad\text{for every }x\in\ell^2(\mathbb Z),
\end{equation}
for almost every $\xi\in[0,1)$.

\subsection{Structure of the matrix $L_{\xi}$}
Fix $\xi\in[0,1)$ and $m_k\in M$. Consider the equation in the unknowns $t\in[0,1)$ and $n\in\mathbb Z$:
\begin{equation}\label{eq:mod}
    \xi + n = t+m_k.
\end{equation}
Equation (4) has the unique solution $(t,n)=(t_k,b_k)$. Since the sequence $(b_k)$ is non-decreasing, there are three regimes in the matrix $L_{\xi}$:
$$\text{I} = \begin{pmatrix}
    \ldots & 0 & \star & \star & 0 & \ldots \\
    \ldots & 0 & \star & \star & 0 & \ldots \\
\end{pmatrix}, \qquad 
\text{II}  = \begin{pmatrix}
    \ldots & 0 & \star & \star & 0 & 0 &\ldots \\
    \ldots & 0 & 0 & \star & \star  & 0 & \ldots \\
\end{pmatrix}, $$
$$
\text{III}  = \begin{pmatrix}
    \ldots & 0 & \star & \star & 0 & 0& 0 & \ldots \\
    \ldots & 0 & 0 & 0 & \star  & \star & 0 &\ldots \\
\end{pmatrix}.$$
The three regimes describe the relative positions of consecutive rows: regime I corresponds to
$b_{k+1}-b_k=0$, regime II to $b_{k+1}-b_k=1$, and regime III to $b_{k+1}-b_k\ge2$. The displayed matrix for regime III represents the case $b_{k+1}-b_k=2$.

\medskip

If there are two neighboring positions with the first regime, then we may isolate a finite segment of rows in such a way that it contains at least as many equations as variables.
$$\newcommand*{\tempb}{\multicolumn{1}{|c}{0}}
\newcommand*{\tempv}{\multicolumn{1}{|c}{\vdots}}
\begin{pmatrix} 
\star & \star &  \tempb &  0  & 0&\ldots  &0\\
\star & \star & \tempb  & 0 & 0 & \ldots  & 0\\
\vdots & \vdots&  \tempv& \vdots  &  \vdots  & & \vdots \\
\star & \star &  \tempb & 0 &  0 & \ldots  &0 \\
\cline{1-2}
0 & \star & \star & 0 & 0 & \ldots  &0 \\
0 & 0 & \star & \star & 0 & \ldots  &0\\
\vdots & &   & \ddots & \ddots& &\vdots  \\
0 & \ldots  & 0 & 0 & \star &  \star  & 0\\
0 & \ldots  & 0 & 0 & 0 & \star & \star  \\
\end{pmatrix}$$
If there are two neighboring places with the third regime, which are not interrupted by the first one, then we may isolate a finite submatrix with strictly more columns than rows:
$$\begin{pmatrix}
\star & \star & 0 & 0 & 0 & 0 & \ldots & \ldots & 0\\
0 & 0 & \star & \star & 0 & 0 &   \ldots & \ldots & 0\\
0 & 0 & 0 & \star & \star & 0 &  \ldots & \ldots & 0\\
\vdots &  &   &   & \ddots & \ddots &  &   &  \vdots \\
0 & \ldots & \ldots & \ldots & 0 & \star & \star & 0 & 0 \\
0 & \ldots & \ldots & \ldots & 0 & 0 & 0 & \star & \star
\end{pmatrix}$$
In particular, this means that $\Ker L_{\xi} \ne \{0\}$. Informally, this means that the first regime helps the frame inequality, while the third regime works against it. The rows of the second regime are the ordinary ones.\par 
We study these two ideas in detail in the next sections.

\section{The operator $L_{\xi}$}
Choose $\xi \in [0,1)$. For every $k \in \Z$ let $t_k \in [0,1)$ be the unique solution of the equation $t+m_k = \xi \bmod 1$  and let
$$b_k = t_k+m_k - \xi  \in \Z.$$
Thus the $k$-th row of $L_{\xi}$ has its non-zero entries in the columns $b_k$ and $b_k+1$. Throughout this section $\delta_0 \in (0,1)$ denotes a separation constant for $M$: now $|x-y|\ge \delta_0$ for any two distinct $x,y \in M$. Whenever condition (2) is assumed with $\alpha=1$, choose $\varepsilon>0$ so that $\sup_k(m_{k+1}-m_k)\leq2-\varepsilon$. When condition (3) is also assumed, one may choose the same $\varepsilon$ as in the condition (3) of Theorem 1.
\subsection{The coefficients}

Fix $n\in\mathbb Z$ and denote by $b_0(n)$ and $b_1(n)$ the coefficients in columns $b_n$ and $b_n+1$ of the $n$-th row, respectively. By the definition of $L_{\xi}$ we have
$$b_0(n) = \begin{cases}
\sinh 2\pi w(\xi-\{m_n\}), & \{m_n\} \leq \xi \\
\sinh 2 \pi w(\xi - \{m_n\}+1 ), & \{m_n\} > \xi
\end{cases} = \sinh 2\pi w t_n,$$
$$b_1(n) = \begin{cases}
\sinh 2\pi w(\xi-\{m_n\}-1), & \{m_n\} \leq \xi \\
\sinh 2 \pi w(\xi - \{m_n\}), & \{m_n\} > \xi
\end{cases} = \sinh 2\pi w (t_n-1).$$

Note that  the number $b_n$ is the smallest integer not less than $m_n-\xi$, since $t_n = b_n + \xi - m_n \in [0,1)$. Hence, $b_n = \lceil m_n - \xi \rceil$. In particular, the sequence $\{b_n\}_{n\in\Z}$ is non-decreasing. We now have
\begin{equation}\label{eq:shift}
b_{n+1}-b_n = \lceil (m_{n+1}-m_n) - t_n  \rceil, \qquad t_{n+1} = t_n + (b_{n+1}-b_n) - (m_{n+1}-m_n).
\end{equation}

\subsection{The blocks}\label{sect:blocks}
In $L_{\xi}$ we may have blocks: pairs of rows of the first regime
$$\begin{pmatrix} 
\ldots & 0 & \star & \star & 0  & \ldots \\
\ldots & 0 & \star & \star & 0 & \ldots 
\end{pmatrix}.$$
Suppose that such rows correspond to $m_j, \ldots, m_{j+k} \in M$. This means that $b_j = \ldots = b_{j+k}$, that is, all these points lie in one interval $(\xi+b-1, \xi+b]$, where $b \in \Z$ is the common value. So we have
\begin{lemma}\label{lem:blocks} For a fixed $\xi \in [0,1)$, the rows corresponding to $m_{j}, \ldots, m_{j+k} \in M$ form a block in $L_{\xi}$ if and only if there exists $n\in \Z$ such that $m_{j}, \ldots, m_{j+k}  \in (\xi + n, \xi + n +1].$
\end{lemma}
\begin{corollary}\label{cor:noblock}
For a fixed $\xi \in [0,1)$ the rows $n$ and $n+1$ form a block if and only if $t_n \ge m_{n+1}-m_n$. They do not form a block if and only if $\xi \in [m_n, m_{n+1})_{\T}$. 
\end{corollary}
\begin{proof}
The rows form a block if and only if $b_{n+1}=b_n$, that is, $  \lceil (m_{n+1}-m_n)-t_n  \rceil  \leq 0$ by \eqref{eq:shift}. Hence, $m_{n+1}-m_n \leq t_n$. The inequality $t_n < m_{n+1}-m_n$ means that $\xi - m_n \bmod 1$ belongs to $[0, m_{n+1}-m_n)$, so $\xi \in [m_n, m_{n+1})_{\T}$.
\end{proof}

Next we show that the blocks are bounded below uniformly with respect to $\xi$.  Without loss of generality we assume that each block has a $2\times 2$ structure, since deleting rows only decreases the left-hand side of \eqref{eq:opercrit}.

\begin{lemma}\label{lem:blockdet}
Let $M \subset \R$ be uniformly discrete with a separation constant $\delta_0 \in (0,1)$. Define $f(\delta) := {{1} \over {2}} \left(\ch 2\pi w (1+\delta) - \ch 2\pi w(1-\delta) \right)$. Then 
\begin{equation}\label{det:block}
    |\det(\text{block})| \ge f(\delta_0) > 0 \quad \text{for any } 2\times 2 \text{ block in } L_{\xi} \text{ and any } \xi \in [0,1).
\end{equation}
 Moreover, if $B$ is a $2\times 2$ block and $B(x_c, x_{c+1})^{T} = (y_n,y_{n+1})^{T}$, then
\begin{equation}\label{eq:blockinv}
\|(x_c,x_{c+1})\| \leq C_1 \|(y_n,y_{n+1})\|, \text{ where } C_1 = {{2\sinh 2\pi w} \over {f(\delta_0)}}.
\end{equation}
\end{lemma}
\begin{proof}
Fix $\xi \in [0,1)$ and consider a block formed by the rows corresponding to $m_n, m_{n+1} \in M$. Put $\delta =  t_n - t_{n+1}$. Since $b_{n+1}=b_n$, the second identity in \eqref{eq:shift} gives $\delta = m_{n+1}-m_n$. We have $\delta \in [\delta_0, 1)$ for any block and any $\xi$ by definition of $\delta_0$. Now we obtain
$$\det \begin{pmatrix}
    b_0(n) & b_1(n) \\ b_0(n+1) & b_1(n+1)
\end{pmatrix} = \sinh 2\pi w t_n \cdot \sinh 2\pi w (t_{n+1}-1) - \sinh 2\pi w (t_n-1) \cdot \sinh 2\pi w t_{n+1} = $$
$$={{1} \over {2}} \left(\ch 2\pi w (\delta-1) - \ch 2\pi w(\delta+1) \right) = -f(\delta),$$
because $\ch$ is an even function. Furthermore, $f(0)=0$ and
$$f'(\delta) =  \pi w \left(\sinh 2\pi w(1+\delta) + \sinh 2\pi w(1-\delta)\right)>0 \qquad \text{for } 0 <\delta < 1.$$
Hence, $f$ is strictly increasing on $[0,1)$ and $|\det(\text{block})| = f(\delta) \ge f(\delta_0)>0$. Neither $\delta$ nor $f(\delta_0)$ depends on $\xi$, which gives the uniformity. Finally, fix a block $B$. All the entries of $B$ are bounded by $\sinh 2\pi w$ in absolute value, so $\|B\| \leq  2 \sinh 2\pi w$. For a $2\times 2$ matrix one has $\|B^{-1}\| = \|B\|/|\det B|$, and \eqref{eq:blockinv} follows.
\end{proof}

\subsection{Cuts} 
In $L_{\xi}$ one may have {\it cuts}, i.e., pairs of rows of the third regime
$$\begin{pmatrix} 
\ldots & 0 & \star & \star & 0 & 0 & 0 & \ldots \\
\ldots & 0 & 0 & 0 & \star & \star & 0 & \ldots 
\end{pmatrix}$$
Suppose such rows correspond to $m_n, m_{n+1}\in M$. By a {\it cut} we mean that the rows $n$ and $n+1$ share no column, that is, $b_{n+1} \ge b_n+2$. Equivalently, there exist $k \in \Z$, $p \ge 2$ and $t_n, t_{n+1}\in[0,1)$ such that $t_n+m_n = \xi+k$ and $t_{n+1}+m_{n+1} = \xi+k+p$. As Corollary \ref{cor:anchor} below shows, $p=2$ whenever condition (2) of Theorem \ref{th:main} holds. So we have

\begin{lemma}\label{lem:div} For a fixed $\xi \in [0,1)$ the $n$-th and $(n+1)$-th rows form a cut in $L_{\xi}$ if and only if 
$$m_{n+1} - m_n>1 \quad \text{ and } \quad \xi \in [m_n, m_{n+1}-1)_{\T}.$$
\end{lemma}

Every cut splits $L_{\xi}$ into the direct sum of two operators acting on disjoint sets of columns. Hence, cuts inside $L_{\xi}$ induce a decomposition 
\begin{equation}\label{eq:segments}
L_{\xi} = \bigoplus_n A_n
\end{equation}
We call such matrices $A_n$ {\it super-segments}. \par
The next corollary shows that the last column of each super-segment and the first column of the next one always contain an entry of modulus at least $\sinh 2\pi w \eps$.

\begin{corollary}\label{cor:anchor}
Let $M$ satisfy condition (2) of Theorem \ref{th:main} with $\alpha=1$ and let the rows $n$ and $n+1$ form a cut in $L_{\xi}$. Then
$$b_{n+1}-b_n = 2, \qquad t_n < 1-\eps, \qquad t_{n+1} \ge \eps.$$
In particular, $|b_1(n)| \ge \sinh 2\pi w \eps$ and $|b_0(n+1)| \ge \sinh (2\pi w \eps)$.
\end{corollary}
\begin{proof}
Put $d_n = m_{n+1}-m_n \le 2-\eps$. By Lemma \ref{lem:div} we have $\xi \in [m_n, m_{n+1}-1)_{\T}$. Hence, $t_n < d_n - 1 \le 1-\eps$ and $|b_1(n)| = \sinh 2\pi w(1-t_n) \ge \sinh 2\pi w \eps$. Moreover, $1 < d_n - t_n \le d_n \le 2-\eps$, so $b_{n+1}-b_n = \lceil d_n - t_n \rceil = 2$ by \eqref{eq:shift}. Therefore $t_{n+1} = t_n + 2 - d_n \ge 2 - d_n \ge \eps$ and $|b_0(n+1)| = \sinh (2\pi w\, t_{n+1}) \ge \sinh (2\pi w \eps)$.
\end{proof}

\subsection{Segments}\label{segments}
Each super-segment $A_n$ from \eqref{eq:segments} consists of blocks (regime I) and 2-diagonal submatrices (regime II) between them. We call the latter {\it $2$-diagonal subsegments}. By {\it a segment} we mean a block together with the $2$-diagonal subsegment following it, up to the next block. Note that there are no rows of the third regime inside $A_n$ by the definition of a super-segment. Thus each $A_n$ splits into segments of the form

$$
\newcommand*{\tempb}{\multicolumn{1}{|c}{0}}
\newcommand*{\tempv}{\multicolumn{1}{|c}{\vdots}}
\begin{pmatrix} 
\star & \star &   \tempb & 0  & 0&\ldots  &0\\
\star & \star &  \tempb  & 0 & 0 &\ldots  & 0\\
\vdots & \vdots& \tempv&  \vdots  &  \vdots  & & \vdots \\
\star & \star & \tempb &  0 &  0 & \ldots  &0 \\
\cline{1-2}
0 & \star & \star & 0 & 0 & \ldots  &0 \\
0 & 0 & \star & \star & 0 & \ldots  &0\\
\vdots & &   & \ddots & \ddots& &\vdots  \\
0 & \ldots  & 0 & 0 & \star &  \star  & 0\\
0 & \ldots  & 0 & 0 & 0 & \star & \star  \\
\end{pmatrix}.$$

Consider a $2$-diagonal subsegment $S$ between two blocks:
$$S = \begin{pmatrix} 
b_0(k) & b_1(k) &   0  &  0 & \ldots & \ldots & 0 \\
0 & b_0(k+1) & b_1(k+1) & 0 & 0 & \ldots  &0 \\
\vdots & &   & \ddots & \ddots& &\vdots  \\
0 & \ldots  & 0 & 0 & b_0(k+n-1) &  b_1(k+n-1)  & 0\\
0 & \ldots  & 0 & 0 & 0 & b_0(k+n) & b_1(k+n) \\
\end{pmatrix}.$$
Recall that $b_0(j) = \sinh 2\pi w t_j$ and $b_1(j) = \sinh 2\pi w (t_j-1)$, where $t_j = \xi - m_j \bmod 1$. 

\begin{lemma}\label{lem:2diag}
    Let the rows $j$ and $j+1$ belong to one $2$-diagonal subsegment of $L_{\xi}$. Then
    $$t_{j+1} - t_j = 1 - (m_{j+1}-m_j).$$
    In particular, $t_{j+1}\ge t_j$ if and only if $m_{j+1}-m_j \leq 1$.
\end{lemma}
\begin{proof}
    The rows $j$ and $j+1$ are of the second regime, so $b_{j+1}=b_j+1$. Hence, \eqref{eq:shift} implies $t_{j+1} = t_j + 1 - (m_{j+1}-m_j)$.
\end{proof}

Fix $\eps>0$ from Theorem \ref{th:main} and put $\eps_0 := \min(\eps/3, \delta_0/2, 1/2)$. By Subsection \ref{sect:blocks}, the inequality $t_j \ge 1-\eps_0$ means that $|b_1(j)| \leq \sinh (2\pi w \eps_0)$, and $t_j \leq \eps_0$ means that $|b_0(j)| \leq \sinh (2\pi w \eps_0)$.

\begin{lemma}\label{lem:order}
Let $M$ satisfy conditions (2) and (3) of Theorem \ref{th:main} with $\alpha=1$. Fix $\xi \in [0,1)$ and let the rows $j_0, j_0+1, \ldots, j_1$ of $L_{\xi}$ form a $2$-diagonal subsegment. Then there exists $r$ with $j_0-1 \leq r \leq j_1$ such that
$$t_j < 1-\eps_0 \quad \text{for } j_0 \leq j \leq r \quad\text{and} \quad t_j > \eps_0 \quad \text{for } r< j \leq j_1.$$
\end{lemma}
\begin{proof}
Let $r$ be the largest index in $[j_0,j_1]$ with $t_r \leq \eps_0$, and put $r=j_0-1$ if there is no such index. The second inequality holds by the choice of $r$. Assume that $t_{j} \ge 1-\eps_0$ for some $j \leq r$. Since $\eps_0<1/2$, we have $t_r \leq \eps_0 < 1-\eps_0$. Hence, $j<r$. Write $d_i = m_{i+1}-m_i$.

Let $n$ be the smallest index in $[j, r-1]$ with $d_n>1$. Observe that such $n$ exists: otherwise $d_i \leq 1$ for all $i \in [j,r-1]$. Hence, $t_{i+1}\ge t_i$ by Lemma \ref{lem:2diag}, and $t_r \ge t_j \ge 1-\eps_0 > \eps_0$, a contradiction. By the same argument applied on $[j,n-1]$ we get $t_n \ge t_j \ge 1-\eps_0$.

Let $k$ be the largest index in $[n, r-1]$ with $d_k>1$. Then $d_i \leq 1$ for $i \in [k+1, r-1]$. Hence, $t_{i+1}\ge t_i$ there. Now $t_{k+1} \leq t_r \leq \eps_0$. Thus we have $n \leq k$, while $\min(d_n, d_k)>1$, $t_n \ge 1-\eps_0$ and $t_{k+1}\leq \eps_0$.

If $n=k$, then the rows $n$ and $n+1$ lie in one $2$-diagonal subsegment, and Lemma \ref{lem:2diag} gives
$$d_n = 1 + t_n - t_{n+1} \ge 1 + (1-\eps_0) - \eps_0 = 2-2\eps_0 > 2-\eps,$$
which contradicts condition (2).

If $n<k$, we show that $\xi$ belongs to the intersection in condition (3). First note that the condition $t_n \ge 1-\eps_0$ means that $\xi \in [m_n - \eps_0, m_n)_{\T}$. Since $d_n>1$ and $\eps_0<\eps$, this set is contained in $[m_n-\eps, m_{n+1}-1)_{\T}$. Second, the condition $t_{k+1}\leq \eps_0$ means that $\xi \in [m_{k+1}-1, m_{k+1}-1+\eps_0]_{\T}$. Since $d_k>1$ and $\eps_0<\eps$, this set is contained in $[m_k, m_{k+1}-1+\eps)_{\T}$. Finally, for every $i$ with $n<i<k$ the rows $i$ and $i+1$ belong to the same $2$-diagonal subsegment and therefore do not form a block. Hence, $\xi \in [m_i,m_{i+1})_{\T}$ by Corollary \ref{cor:noblock}.
\end{proof}

\section{Necessity in Theorem \ref{th:main}}

\subsection{Necessity of (1)} The proof follows the same lines as in the classical case of Bernstein spaces. In the Gabor analysis setting, one can apply Theorem 1.5 from \cite{AFK} since the window function $g$ belongs to Feichtinger's algebra $M^1(\mathbb{R})$. Indeed, let $\varphi(t)=e^{-\pi t^2}$ and let $V_{\varphi}g(x,\omega) = \int_{\R} g(t)\varphi(t-x)e^{-2\pi i \omega t}\,dt$ be the short-time Fourier transform of $g$. Put $h_x(t) := g(t)\varphi(t-x)$. By the Cauchy--Schwarz inequality and Plancherel's theorem we have
$$\int_{\R}|V_{\varphi}g(x,\omega)|\,d\omega \leq \left(\int_{\R}{{d\omega}\over{1+\omega^2}}\right)^{1/2} \left(\int_{\R}(1+\omega^2)|\widehat{h_x}(\omega)|^2 d\omega \right)^{1/2} \lesssim \|h_x\|_2 + \|h_x'\|_2 .$$
Since $|g(t)|+|g'(t)| \lesssim (1+t^2)^{-1}$ and both $\varphi$ and $\varphi'$ decay faster than any power, we have $\|h_x\|_2 + \|h_x'\|_2 \lesssim (1+x^2)^{-1}$. Integrating in $x$ we obtain $\|V_{\varphi}g\|_{L^1(\R^2)} < \infty$, that is, $g \in M^1(\R)$.

\subsection{Reduction to a uniformly discrete subset}
We first need the stability of irregular Gabor frames under a uniformly small perturbation of the time-frequency nodes.
\begin{lemma}\label{lem:stab}
Let $g \in M^1(\R)$ and let $M \subset \R$ be such that $\G(g, M \times \Z)$ is a frame in $L^2(\R)$. Then there exists $\delta_1>0$ with the following property: for every map $\sigma \colon M \to \R$ with $\sup_{m \in M}|\sigma(m)-m| \leq \delta_1$ the system $\G(g, \sigma(M)\times \Z)$ is a frame in $L^2(\R)$.
\end{lemma}
\begin{proof}
Put $\Lambda := M \times \Z$ and assume the contrary. Then for every $n \in \N$ there exists $\sigma_n \colon M \to \R$ with $\sup_{m \in M}|\sigma_n(m)-m| \leq 1/n$ such that $\G(g, \Lambda_n)$ is not a frame, where $\Lambda_n = \tau_n(\Lambda) = \sigma_n(M)\times\Z$ for $\tau_n(m,k) = (\sigma_n(m), k)$. Since $\G(g,\Lambda)$ is a frame, the set $\Lambda$ is relatively separated by \cite[Lemma 2.2]{GOR}. As $\sup_{\lambda \in \Lambda}|\tau_n(\lambda)-\lambda| \leq 1/n \to 0$, the sets $\Lambda_n$ are Lipschitz convergent to $\Lambda$ by \cite[Example 6.2]{GOR}, and therefore $\G(g,\Lambda_n)$ is a frame for all sufficiently large $n$ by \cite[Theorem 7.1]{GOR}.
\end{proof}

Lemma \ref{lem:stab} yields the standard reduction to a uniformly discrete subset. Fix $\delta \in (0,\delta_1)$ and choose a maximal $\delta$-separated subset $N \subseteq M$ by Zorn's lemma. Now for every $m \in M$ there is $\phi(m) \in N$ with $|m-\phi(m)|<\delta$, since $N$ is maximal. Put $N' := \phi(M) \subseteq N$. The map $\phi \colon M \to \R$ moves every point by less than $\delta \leq \delta_1$ and $\phi(M) = N'$, so $\G(g, N'\times \Z)$ is a frame by Lemma \ref{lem:stab}. Since $N' \subseteq N \subseteq M$, the system $\G(g, N \times \Z)$ is a frame as well.\par

Moreover, since $\G(g, N \times \Z)$ is a frame and $g \in M^1(\R)$, the planar lower Beurling density of $N \times \Z$ is greater than $1$ by \cite[Theorem 1.5]{AFK}. Hence condition (1) holds for $N$.\par
Thus we may assume from now on that $M$ is uniformly discrete with a separation constant $\delta_0 \in (0,1)$. Let $c>0$ denote a constant for which \eqref{eq:opercrit} holds for the system $\G(g, M \times \Z)$ and let $\delta_1>0$ denote the constant of Lemma \ref{lem:stab} for this system. It remains to prove condition (2) and to find $\eps>0$ for which condition (3) holds.

\subsection{Necessity of (2)} 

First assume that $m_{k+1}-m_k > 2$ for some $k \in \Z$. Then for any $\xi \in (m_k+1, m_{k+1}-1)_{\T} \ne \emptyset$ we have an identically zero column in the matrix of the operator $L_{\xi}$, which contradicts inequality \eqref{eq:opercrit}. Hence $d_k = m_{k+1}-m_k \leq 2$ for all $k \in \Z$.\par
Now fix $k$ with $d_k \in (1,2]$. Without loss of generality suppose $m_k = 0$ and $m_{k+1} = d_k$. Fix $\xi \in (0, d_k-1)$. By Lemma \ref{lem:div}, the rows $k$ and $k+1$ form a cut. Moreover, $\xi - d_k \in (-2,-1)$, so
$$t_{k+1} = \xi - d_k \bmod 1 = 2 - d_k + \xi.$$

Let $l = b_{k+1}$ and let $e_l \in \ell^2(\Z)$ be the corresponding coordinate vector. Since the rows $k$ and $k+1$ form a cut, no row has its second entry in the $l$-th column. The rows $j$ with $b_j = l$ are consecutive, start with the row $k+1$ and satisfy $t_{j+1} = t_j - d_j$ by \eqref{eq:shift}. Hence $t_j \leq t_{k+1}$ for all of them and their number is at most $1 + 1/\delta_0$. Now for almost every $\xi \in (0, d_k-1)$ we have
$$c = c\cdot \|e_l\|^2 \leq \|L_{\xi} e_l\|^2 = \sum_{j\,:\, b_j = l} \sinh^2 2\pi w\, t_j \leq \left(1+{{1}\over{\delta_0}}\right) \sinh^2 2\pi w (2-d_k+\xi)$$
by \eqref{eq:opercrit}. Now by taking $\xi$ arbitrarily small we obtain
$$\sinh 2\pi w (2-d_k) \ge \sqrt{{{c\,\delta_0}\over{1+\delta_0}}}.$$
Hence there exists $\eps_2 \in (0,1]$, depending only on $c$, $w$ and $\delta_0$, such that $m_{k+1}-m_k \leq 2-\eps_2$ for every $k \in \Z$. In particular, $\sup_k |m_{k+1} - m_k| < 2$, which is condition (2).\par

\subsection{Necessity of (3)}

\begin{lemma}\label{lem:auxdiv} Fix $\xi \in [0,1)$ and suppose that in $L_{\xi}$ the rows $k, k+1$ and the rows $j, j+1$ form cuts for some $k<j$. Assume that no rows $i, i+1$ with $k<i<j$ form a block. Then $\Ker L_{\xi} \ne \{0\}$.
\end{lemma}
\begin{proof}
Consider the rows $k+1, \ldots, j$. Their non-zero entries lie in the columns $b_{k+1}, \ldots, b_j+1$. Since no rows $i, i+1$ with $k<i<j$ form a block, we have $b_{i+1} \ge b_i+1$ for these $i$, whence
$$b_j+1 - b_{k+1} + 1 \ge (j-k-1)+2 > j-k.$$
Thus the number of columns strictly exceeds the number of rows, and there exists $x \ne 0$ supported on the columns $b_{k+1}, \ldots, b_j+1$ such that the rows $k+1, \ldots, j$ of $L_{\xi}x$ vanish. \par 
The remaining rows of $L_{\xi}x$ vanish automatically. Indeed, the rows $i \leq k$ have their non-zero entries in the columns $b_i, b_i+1$ with $b_i + 1 \leq b_k+1 \leq b_{k+1}-1$, and the rows $i \ge j+1$ have their non-zero entries in the columns with indices at least $b_{j+1} \ge b_j + 2$. Hence $L_{\xi}x=0$.
\end{proof}

Put $\eps_3 := \min(\delta_0/2,\delta_1)$. Let $n<k$ be indices such that $m_{n+1}-m_n>1$ and $m_{k+1}-m_k>1$, and let
$$\Omega_{\eps_3} = [m_n-\eps_3,\, m_{n+1}-1)_{\T} \cap \left( \bigcap_{j=n+1}^{k-1} [m_j, m_{j+1})_{\T} \right) \cap [m_k,\, m_{k+1}-1+\eps_3)_{\T}.$$
Assume for contradiction that $\Omega_{\eps_3} \ne \emptyset$. Put
$$M' := \left(M \setminus \{m_n, m_{k+1}\}\right) \cup \{m_n-\eps_3,\ m_{k+1}+\eps_3\},$$
so that $m_n' = m_n - \eps_3$, $m_{k+1}' = m_{k+1}+\eps_3$ and $m_j'=m_j$ for $j \ne n, k+1$. The points of $M'$ are strictly increasing in $j$, and $M'$ is uniformly discrete. We have $\sup_j|m_j'-m_j| = \eps_3 \leq \delta_1$, so $\G(g, M'\times \Z)$ is a frame by Lemma \ref{lem:stab}. We obtain
$$\Omega_{\eps_3} = [m_n',\, m_{n+1}'-1)_{\T} \cap \left( \bigcap_{j=n+1}^{k-1} [m_j', m_{j+1}')_{\T} \right) \cap [m_k',\, m_{k+1}'-1)_{\T}.$$
Fix $\xi \in \Omega_{\eps_3}$ and consider the operator $L_{\xi}$ associated with $M'$. Since $\xi \in [m_n', m_{n+1}'-1)_{\T}$ and $\xi \in [m_k', m_{k+1}'-1)_{\T}$, the rows $n, n+1$ and the rows $k, k+1$ form cuts by Lemma \ref{lem:div}. For every index $n<j<k$ we have $\xi \in [m_j', m_{j+1}')_{\T}$, so the rows $j$ and $j+1$ do not form a block by Corollary \ref{cor:noblock}. Hence these two cuts are not interrupted by a block and $\Ker L_{\xi}\ne\{0\}$ by Lemma \ref{lem:auxdiv}.\par
Finally, $\Omega_{\eps_3}$ is a non-empty intersection of finitely many arcs of the form $[a,b)_{\T}$, hence $|\Omega_{\eps_3}|>0$, which leads us to the contradiction.\par

It remains to note that conditions (2) and (3) become weaker as $\eps$ decreases, so it suffices to verify (2) for some $\eps_2>0$ and (3) for some $\eps_3>0$: they then both hold with $\eps := \min(\eps_2,\eps_3)$.

\section{Sufficiency}

Since $D_+(M)<\infty$ and $g \in M^1(\R)$, it is sufficient to consider only a uniformly discrete set $M$ with a separation constant $\delta_0 \in (0,1)$ such that (1)--(3) hold with some $\eps>0$. We keep $\eps_0 = \min(\eps/3,\delta_0/2)$ from Section \ref{segments}. Fix $\xi \in [0,1)$ and consider the operator $L_{\xi}$. \par 
Removing rows from the blocks only decreases the left-hand side of \eqref{eq:opercrit}, so we may assume that every block is $2\times2$. In each maximal block, retain its first and last rows
and delete all intermediate rows. If a deleted point lies between two rows which are consecutive after the deletion, then by Lemma \ref{lem:blocks} these two rows would belong to one block. Hence every cut pair and every pair of consecutive rows of a $2$-diagonal subsegment of the reduced operator consists of consecutive points of the original set $M$, and conditions (2)--(3), Corollary \ref{cor:anchor} and Lemmas \ref{lem:order} and \ref{lem:boundary} below apply to them. \par
Consider a decomposition of $L_{\xi}$ into super-segments $L_{\xi} = \bigoplus_n A_n$ (see \eqref{eq:segments}).  Each $A_n$ is a concatenation of blocks and of the $2$-diagonal subsegments between them (see Subsection \ref{segments}). A subsegment adjacent to a cut has a block on one side only. Since the summands in \eqref{eq:segments} act on pairwise disjoint sets of columns, it suffices to show that the operators $A_n$ are bounded below uniformly in $n$ and $\xi$. \par 

\begin{lemma}\label{lem:runs}
There exists $N_0\in \N$ such that for almost every $\xi \in [0,1)$ every $2$-diagonal subsegment of $L_{\xi}$ contains at most $N_0$ rows.
\end{lemma}
\begin{proof}
Since $D_{-}(M)>1$, one may find $c>0$ and $a_0>0$ such that 
$$\#(M\cap[R,R+a])\ge (1+c)a \quad \text{ for all } R \in \R \text{ and all } a \ge a_0.$$
Let the rows $j, \ldots, j+L$ form a $2$-diagonal subsegment. Summing the identity of Lemma \ref{lem:2diag} we get $m_{j+L}-m_j = L+t_j-t_{j+L}>L-1$, while $m_j, \ldots, m_{j+L}$ are consecutive points of $M$. Hence, $[m_j, m_{j+L}]$ contains exactly $L+1$ points of $M$. \par 
If $L-1\geq a_0$, this gives $L+1\geq(1+c)(L-1)$, hence $L\leq1+2/c$.
Otherwise $L<a_0+1$. Since the subsegment contains $L+1$ rows, we may take
\[
N_0:= \left\lceil \max\left\{ a_0+2,\, 2+\frac{2}{c}\right\} \right\rceil.
\]

\end{proof}

{\bf Step 1}. Every super-segment contains a block. Assume that $A_n$ contains no block. By Lemma \ref{lem:runs} it has at most $N_0$ rows. Hence it is bounded by cuts on both sides: the rows $k, k+1$ and the rows $j, j+1$ form cuts for some $k<j$, and no rows $i, i+1$ with $k<i<j$ form a block. By Lemma \ref{lem:div} we have 
$$m_{k+1}-m_k>1, \quad m_{j+1}-m_j>1 \quad  \text{ and } \quad \xi \in [m_k, m_{k+1}-1)_{\T} \cap [m_j, m_{j+1}-1)_{\T},$$
while $\xi \in [m_i,m_{i+1})_{\T}$ for $k<i<j$ by Corollary \ref{cor:noblock}. Hence $\xi$ lies in the intersection of condition (3) for the pair $k<j$, which contradicts (3).

{\bf Step 2}. Let the rows $j_0,\ldots,j_1$ form a $2$-diagonal subsegment. By Lemma \ref{lem:order} there is $r \in [j_0-1,j_1]$ with
$$|b_1(j)| = \sinh 2\pi w(1-t_j) \ge \sinh 2\pi w \eps_0 \ \ (j_0 \leq j \leq r), \qquad |b_0(j)| = \sinh 2\pi w\,t_j \ge \sinh 2\pi w \eps_0 \ \ (r<j\leq j_1).$$
Thus along a subsegment we move from a large $|b_1|$ to a large $|b_0|$. If a subsegment is adjacent to a cut, this dichotomy degenerates.

\begin{lemma}\label{lem:boundary}
Let $M$ satisfy conditions (2) and (3) of Theorem \ref{th:main} with $\alpha=1$, fix $\xi \in [0,1)$ and let the rows $j_0, \ldots, j_1$ of $L_{\xi}$ form a $2$-diagonal subsegment.
\begin{enumerate}
\item If the rows $j_0-1$ and $j_0$ form a cut, then $t_j > \eps_0$ for all $j_0 \leq j \leq j_1$.
\item If the rows $j_1$ and $j_1+1$ form a cut, then $t_j < 1-\eps_0$ for all $j_0 \leq j \leq j_1$.
\end{enumerate}
\end{lemma}
\begin{proof}
We prove the first item, the proof of the second item is symmetric. Write $d_i = m_{i+1}-m_i$ and put $n := j_0-1$. By Lemma \ref{lem:div} we have $d_n>1$ and $\xi \in [m_n, m_{n+1}-1)_{\T}$, and by Corollary \ref{cor:anchor} we have $t_{j_0} \ge \eps > \eps_0$.\par
Assume that $t_j \leq \eps_0$ for some $j \in [j_0,j_1]$; then $j>j_0$. Let $k$ be the largest index in $[j_0, j-1]$ with $d_k>1$. Such $k$ exists: otherwise $d_i \leq 1$ for all $i \in [j_0,j-1]$, hence $t_{i+1}\ge t_i$ by Lemma \ref{lem:2diag} and $t_j \ge t_{j_0}>\eps_0$. By maximality $d_i \leq 1$ for $i \in [k+1,j-1]$, whence $t_{k+1}\leq t_j \leq \eps_0$.\par
Thus $n<k$, $d_n>1$ and $d_k>1$. Furthermore, $\xi \in [m_n, m_{n+1}-1)_{\T} \subset [m_n-\eps, m_{n+1}-1)_{\T}$; the inequality $t_{k+1}\leq \eps_0$ means that $\xi \in [m_{k+1}-1, m_{k+1}-1+\eps_0]_{\T} \subset [m_k, m_{k+1}-1+\eps)_{\T}$, since $d_k>1$ and $\eps_0<\eps$; and for every $i$ with $n<i<k$ the rows $i$ and $i+1$ lie in one $2$-diagonal subsegment, so they do not form a block and $\xi \in [m_i, m_{i+1})_{\T}$ by Corollary \ref{cor:noblock}. Hence the intersection in condition (3) for the pair $n<k$ is non-empty, a contradiction.
\end{proof}

{\bf Step 3}. Fix $n \in \Z$ and consider the equation $A_nx = y$. The super-segment $A_n$ consists of $2\times2$ blocks $\{B_k\}_k$ and of $2$-diagonal subsegments $\{S_k\}_k$ between them. Assign to every block its two columns and to every subsegment those of its columns which are not shared with a neighbouring block. Then every row of $A_n$ belongs to exactly one group, while every column belongs to at least one group. Hence,
$$\|x\|^2 \leq \sum_k \|x_{B_k}\|^2 + \sum_k \|x_{S_k}\|^2, \qquad \|y\|^2 = \sum_k \|y_{B_k}\|^2 + \sum_k \|y_{S_k}\|^2.$$

The rows of a block meet no columns outside it, so Lemma \ref{lem:blockdet} gives
$$\|x_{B_k}\|^2 \leq C_1^2 \|y_{B_k}\|^2,\text{ where }C_1 = {{2\sinh 2\pi w}\over{f(\delta_0)}}.$$
Consider a subsegment $S$ with rows $1, \ldots, q$, where $q \leq N_0$ by Lemma \ref{lem:runs}. Let $x_0, \ldots, x_q$ be the columns met by these rows, so that
\begin{equation}\label{eq:2diag}
b_0(j)x_{j-1} + b_1(j)x_j = y_j, \qquad 1 \leq j \leq q.
\end{equation}
Equations \eqref{eq:2diag} may be resolved in two directions:
$$\begin{cases} x_j = {{y_j}\over{b_1(j)}} - {{b_0(j)}\over{b_1(j)}}x_{j-1}, & \text{if } |b_1(j)| \ge \sinh 2\pi w \eps_0, \\
x_{j-1} = {{y_j}\over{b_0(j)}} - {{b_1(j)}\over{b_0(j)}}x_j, & \text{if } |b_0(j)| \ge \sinh 2\pi w \eps_0.
\end{cases}$$
All the entries of $L_{\xi}$ are bounded by $\sinh 2\pi w$ in absolute value, so each step multiplies the previous unknown by a factor of modulus at most $\sinh 2\pi w/\sinh 2\pi w \eps_0$. Iterating at most $N_0$ times and applying the Cauchy--Schwarz inequality we obtain, for every column $x_j$ reached by such a recursion started at a column $x_{\ast}$,
$$|x_j| \leq {{Q_0 \sqrt{N_0}}\over{\sinh 2\pi w \eps_0}}\|y_S\| + Q_0|x_{\ast}|, \qquad Q_0 := \left({{\sinh 2\pi w}\over{\sinh 2\pi w \eps_0}}\right)^{N_0}.$$
It remains to check that in each case all the columns of $S$ are reached.\par
 If $S$ lies between two blocks, then $x_0$ and $x_q$ belong to these blocks, and by Step 2 the forward recursion started at $x_0$ produces $x_1, \ldots, x_r$, while the backward one started at $x_q$ produces $x_{q-1}, \ldots, x_r$. Together they cover $x_1, \ldots, x_{q-1}$. \par 
If $S$ is preceded by a cut, then $x_0$ is a column of $S$, while $x_q$ belongs to a block (which exists by Step 1). By Lemma \ref{lem:boundary} (item 1) we have $|b_0(j)| \ge \sinh 2\pi w \eps_0$ for all $j$, so the backward recursion started at $x_q$ produces $x_{q-1}, \ldots, x_0$. \par 
The case of a subsegment followed by a cut is symmetric. In all three cases
$$\|x_{S}\|^2 \leq C_2\left(\|y_{S}\|^2 + \|x_{B}\|^2 + \|x_{B'}\|^2\right),$$
where $C_2$ depends only on $w, \eps_0$ and $N_0$, and $B, B'$ are the blocks adjacent to $S$ (one of them being absent if $S$ is adjacent to a cut). Since every block is adjacent to at most two subsegments, we conclude that
$$\|x\|^2 \leq \left(C_1^2 + 2C_2C_1^2\right)\sum_k \|y_{B_k}\|^2 + C_2 \sum_k \|y_{S_k}\|^2 \leq C_3\|y\|^2,$$
where $C_3$ is independent of $n$ and $\xi$. The super-segments act on pairwise disjoint sets of columns, so $\|x\|^2 \leq C_3\|L_{\xi}x\|^2$ for every $x \in \ell^2(\Z)$ and almost every $\xi \in [0,1)$. Hence the proof is complete by Theorem \ref{main:crit}.


\begin{thebibliography}{99}
\bibitem{AFK} G. Ascensi, H. G. Feichtinger, and N. Kaiblinger, {\it Dilation of the Weyl symbol and Balian-Low Theorem}, Trans. of the American Math. Society, Vol. 366, Number 7, July 2014, pp. 3865--3880.

\bibitem{BB}A. Baranov, Y. Belov, {\it Irregular sampling for hyperbolic secant type functions}, Adv. Math. 458, Part B, Article ID 109981, 16 p. (2024).

\bibitem{BKL2}  Yu. Belov, A. Kulikov, Yu. Lyubarskii, \textit{Irregular Gabor frames of Cauchy kernels}, Appl. Comput. Harmon. Anal. 57 (2022), pp. 101--104.


\bibitem{BKL1} Y. Belov, A. Kulikov, Y. Lyubarskii, \textit{Gabor frames for rational functions}, Inventiones Mathematicae, 231:431--466 (2023).

\bibitem{BK} Y. Belov and A. Kulikov, \textit{Gabor frames for functions supported on a semi-axis}, Applied and Computational Harmonic Analysis, volume 85, 101892 (2026).

\bibitem{BAVS} Yu. Belov and A. V. Semenov, \textit{Frame set for shifted sinc-function},  Applied and Computational Harmonic Analysis, 71 (2024), 101654.

\bibitem{PGLST}J. de Dios Pont, K.~Gr\"ochenig, L.~Liehr, I. Shafkulovska, M. Taylor, {\it Gabor frames of totally positive functions: a complete characterization}, https://arxiv.org/abs/2608.04992.


\bibitem{Gro} K.~Gr\"ochenig, {\em  Foundations of Time-Frequency Analysis,} Birkh\"auser, Boston, MA, 2001.

\bibitem{Mystery} K.~Gr\"ochenig, {\em The mystery of Gabor frames},
Journal of Fourier Analysis and Applications, 20 (2014), no. 4, pp. 865--895.

\bibitem{GK} K.~Gr\"ochenig, S.~Koppensteiner, {\em Gabor frames: characterizations and
coarse structure}, in: New Trends in Applied Harmonic Analysis, Vol.~2, Applied and Numerical
Harmonic Analysis, Birkh\"auser/Springer, Cham, 2019, pp. 93--120.

\bibitem{GOR} K. Gr\"ochenig, J. Ortega-Cerd\`a, and J. L. Romero, \textit{Deformation of Gabor systems}, Adv. Math. \textbf{277} (2015), 388--425.

\bibitem{Gro2} K.~Gr\"ochenig, J.~St\"ockler, {\em Gabor frames and totally positive functions,}
Duke Mathematical Journal, 162 (6), 1003--1031, (2013).

\bibitem{Gro1}K.~Gr\"ochenig, J.L.~Romero, J.~St\"ockler, {\em Sampling theorems for shift-invariant spaces, Gabor frames, and totally positive functions,}
Inventiones mathematicae, 211 (3), 1119--1148, (2018).

\bibitem{Heil} C. Heil, {\it History and evolution of the density theorem for Gabor frames},
Journal of Fourier Analysis and Applications, 13 (2007), no. 2, pp. 113--166.


\bibitem{Jans2} A. J. E. M.~Janssen, {\em Some Weyl-Heisenberg frame bound calculations,} Indag. Math., 7:165--183, (1996).

\bibitem{Jans3} A. J. E. M.~Janssen, {\em On generating tight Gabor frames at critical density,} J. Fourier Anal.
Appl., 9(2):175--214, (2003).


\bibitem{JansStr} A.~Janssen, T.~Strohmer, {\em Hyperbolic secants yield Gabor frames},
Appl. Comput. Harmon. Anal., 12, 259--267, (2002).

\bibitem{L} Yu. Lyubarskii, {\em Frames in the Bargmann space of entire functions,} Entire and Subharmonic Functions, Adv. Soviet Math., vol. 11, Amer.\ Math.\ Soc., Providence, RI, 1992, pp. 167--180.

\bibitem{RS} J. Ramanathan and T. Steger, {\it Incompleteness of Sparse Coherent States}, Applied and Computational Harmonic Analysis, Volume 2, Issue 2, April 1995, pp. 148--153.

\bibitem{S} K.~Seip, {\em Density theorems for sampling and interpolation in the Bargmann--Fock space. I,} J. Reine Angew.\ Math.\ {\bf 429} (1992) 91--106.

\bibitem{SW}
K.~Seip and R.~Wallst\'en, {\em Density theorems for sampling and interpolation in the Bargmann--Fock space. II,} J. Reine Angew.\ Math.\ {\bf 429} (1992), 107--113.

\bibitem{AVS} A. V. Semenov, {\it Universal frame set for rational functions}, Applied and Computational Harmonic Analysis, 86 (2027), 101915.
\end{thebibliography}
\end{document}